\documentclass[12pt]{article}
\usepackage[cp1251]{inputenc}
\usepackage{epsfig,amssymb,epic,eepic,color,graphics,amsmath,amsthm}
\usepackage[russian,english]{babel}
\usepackage{multicol}
\usepackage{enumerate}

\newtheorem{theorem}{Theorem}

\newtheorem{lemma}{Lemma}
\newtheorem{proposition}{Proposition}

\begin{document}

\title{
\medskip
\Large On embedding of arcs and circles in 3-manifolds in an application to dynamics of rough 3-diffeomorhisms with two-dimensional expanding attractor\footnote{This work was supported  by the Russian Science Foundation  (project 17-11-01041),  except for the section 5 (Application for dynamics.  Proof of the theorem 3) which was supported  by the Russian Science Foundation  (project  14-41-00044).}}
\author{V.\,Z. Grines, E.\,V. Kruglov, T.\,V. Medvedev, O.\,V. Pochinka}
\date{}

\maketitle

\begin{abstract}
A topological classification of many classes of dynamical systems with regular dynamics in low dimensions is often reduced to combinatorial invariants. In dimension 3 combinatorial invariants are proved to be insufficient even for simplest Morse-Smale diffeomorphisms. The complete topological invariant for the systems with a single saddle point on the 3-sphere is the embedding of the homotopy non-trivial knot into the manifold $\mathbb S^2\times\mathbb S^1$. If a diffeomorphism has several saddle points their unstable separatrices form arcs frames in the basin of the sink and circles frame in the orbits space. Thus, the type of embedding of the circles frame into $\mathbb S^2\times\mathbb S^1$ is a topological invariant for diffeomorphisms of this kind and this type turns out to be the complete topological invariant for some classes of Morse-Smale 3-diffeomorphisms. Recently it was shown that the problem of embedding of a 3-diffeomorphism into a topological flow is interconnected with the properties of embedding of the arcs frame into the 3-Euclidean space. In this paper we consider the criteria for the tame embedding of an arcs frame into the 3-Euclidean space as well as for the trivial embedding of circles frame into $\mathbb S^2\times\mathbb S^1$. We apply this criteria to prove that frames of one-dimensional separatrices in basins of sources of rough 3-diffeomorhisms with two-dimensional expanding attractor are tamely embedded and their spaces of orbits are trivial embeddings of  circles frame into $\mathbb S^2\times\mathbb S^1$.
\end{abstract}

\section{Introduction and formulation of results}

Let $\mathbb R^3$ be 3-dimensional Euclidean space, $\mathbb{S}^{2}=\{(x_1,x_2,x_3)\in \mathbb{R}^3: x^2_1+x^2_2+x^2_3=1\}$ and $O(0,0,0)\in\mathbb R^3$ be the coordinates origin. We represent the manifold $\mathbb R^3\setminus \{O\}$ as $\mathbb S^2\times\mathbb R$ if we assign the point $p(x)=(q(x),r(x))\in \mathbb{S}^{2}\times \mathbb{R}$ to every point $x=(x_1,x_2,x_3)\in \mathbb{R}^3\setminus \{O\}$  where $$q(x)=\left(\frac{x_1}{||x||},\frac{x_2}{||x||},\frac{x_3}{||x||}\right),~~r(x)=\log_{2}||x||,~~||x||=\sqrt{x_1^2+x_2^2+x_3^2}.$$ We say that $\mathbb L\subset(\mathbb R^3\setminus O)$ is an {\it open ray from $O$} if $p(\mathbb L)=\{s\}\times \mathbb{R}$ for some point $s\in\mathbb S^2$. A one-dimensional submanifold $L\subset (\mathbb R^3\setminus O)$ is called an {\it infinite arc from $O$} if $r(L)=\mathbb R$.  

For $\nu\in\mathbb N$ and  a collection $L_1,...,L_\nu$ of pairwise infinite arcs from $O$, the union $${F}_\nu=\bigcup\limits_{j=1}^{\nu} L_j\cup O$$ to be the {\it frame} of $\nu$ arcs. Everywhere below we will assume that every arc of the frame $F_\nu$ is invariant with respect to the homothety $A:\mathbb{R}^3\to \mathbb{R}^3$ defined by $$A(x_1,x_2,x_3)=(2x_1,2x_2,2x_3).$$ A frame  $$\mathbb{F}_\nu=\bigcup\limits_{j=1}^{\nu}\mathbb L_j\cup O$$ of pairwise distinct open rays $\mathbb L_1,...,\mathbb L_\nu$ from $O$ to be {\it standard frame}. We say a frame $F_\nu$ to be {\it tame} if there exists a homeomorphism $H:\mathbb{R}^3\to \mathbb{R}^3$ such that  ${F}_\nu=H(\mathbb F_\nu)$; otherwise we say the frame  $F_\nu$ to be {\it wild}.

\begin{figure}[h]
\centerline{\includegraphics[width=7 true cm]{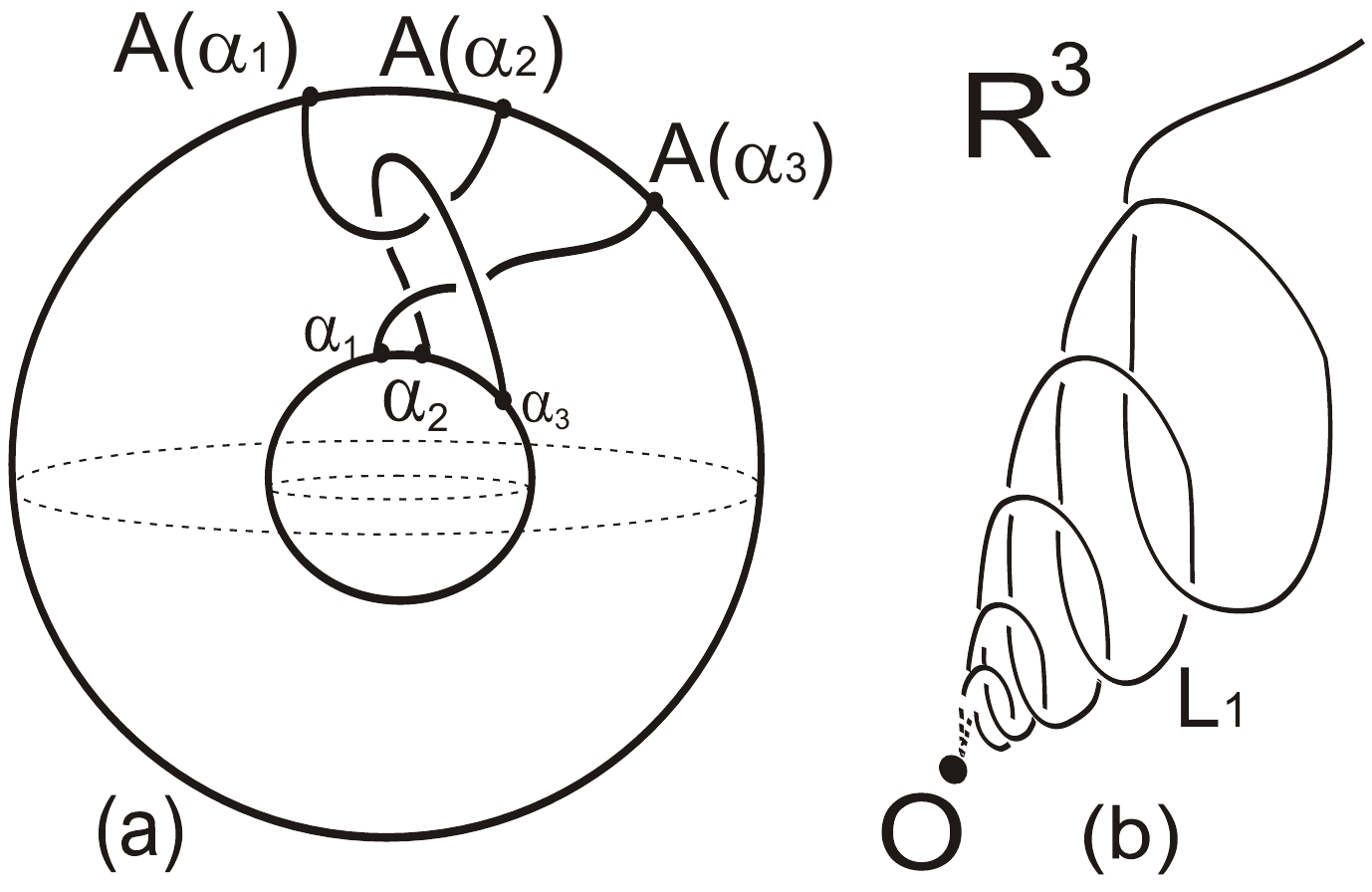}}
\caption{Construction of wild arcs in $\mathbb R^3$}
\label{wildarcconstruction}
\end{figure} 

If $\nu=1$ then the frame is reduced to one arc. E.~Artin and R.~Fox \cite{ArFo} were the first to construct an example of a wild arc in 1948  (see Fig \ref{wildarcconstruction} (b)). Figure \ref{wildarcconstruction} (a) shows the intersection of the arc $L_1$ with the 3-annulus $\{(x_1,x_2,x_3)\in\mathbb R^3:1\leq x_1^2+x_2^2+x_3^2\leq 4\}$. The arc $L_1$ is the union of all positive and negative iterations of all the arcs in this annulus by the homothety $A$.  

In 1955 O.~Harrold, H.~Griffith and E.~Posey E.E. in \cite{HGP1955} proved the following criterion for an arc $F_1=L_1\cup O$ to be tame.

\begin{proposition}[\cite{HGP1955}, Theorem VII]\label{tame-arc} An arc $F_1$ is tame if and only if there is a 3-ball $B(O)$ with $O$ in its interior such that $\partial B(O)\cap L_1$ is a single point.
\end{proposition} 

\begin{figure}[h]
\centerline{\includegraphics[width=11 true cm]{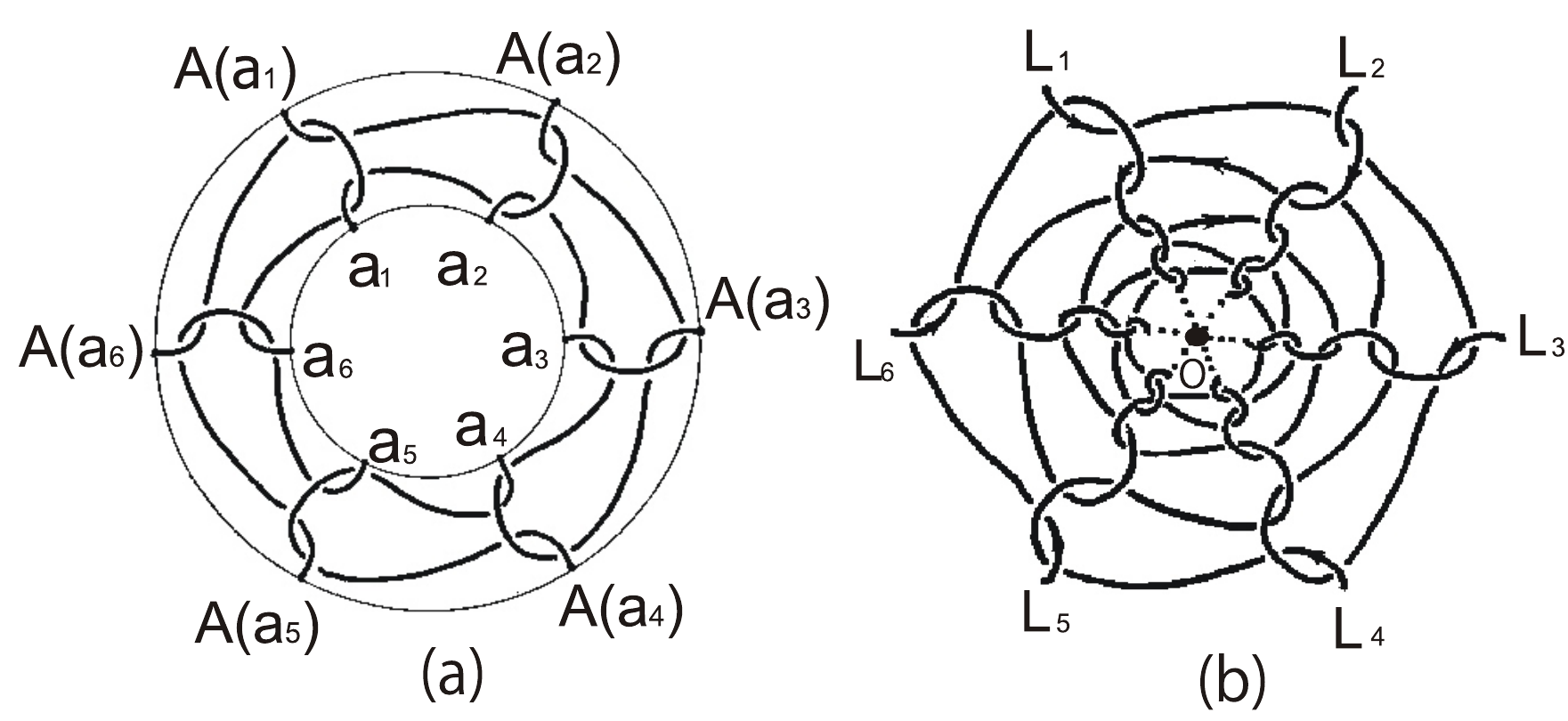}}
\caption{Debrunner-Fox example for $\nu=5$}
\label{mildwild}
\end{figure}

Notice that if each member of the frame $F_\nu\subset \mathbb{R}^3,\nu>1$ is tame this does not necessarily mean that the frame is tame. H.~Debrunner and R.~Fox in 1960 in \cite{DF} constructed an example of a {wild frame} of $\nu>1$ arcs each of which was tame (see Fig. \ref{mildwild}). Moreover, in this frame $F_\nu=\bigcup\limits_{j=1}^\nu L_j\cup O$ 2-sphere $\mathbb S^2$ intersects $L_j$ at exactly one point for each $j\in \{1,\dots,\nu\}$. 

Thus, Harrold-Griffith-Posey criterion cannot be generalized for $\nu>1$. On the other hand, for a frame $F_\nu$ (for any $\nu$) to be tame it is necessary that there exists a foliation $P$ of the set $\mathbb R^3\setminus O$ whose leaves are 2-spheres and every one of these spheres intersects each $L_j$ at a single point. Indeed, from the definition of the tame frame as $F_\nu=H(\mathbb F_\nu)$ it follows that such a foliation is formed by the images by $H$ of the spheres $$\mathbb P_t=\{(x_1,x_2,x_3)\in \mathbb{R}^3: x^2_1+x^2_2+x^2_3=t^2\},\,{t>0}$$ which form the foliation $\mathbb P$.

In this paper we prove the following criterion for an arcs frame $F_\nu$ to be tame (for any $\nu\in\mathbb N$).

\begin{theorem}\label{fibr} A frame $F_\nu=\bigcup\limits_{j=1}^\nu L_j\cup O$ is tame if and only if there is a homeomorphism $G:\mathbb R^3\setminus \{O\}\to\mathbb R^3\setminus \{O\}$ such that every sphere of the foliation $P=G(\mathbb P)$ intersects each arc $L_j,j\in\{1,\dots,\nu\}$ at a single point.
\end{theorem}

For every standard arcs frame $\mathbb F_\nu$ in $\mathbb R^3$ there is the circles frame in $\mathbb S^2\times\mathbb S^1$ corresponding to it by means of the cover $\pi: \mathbb{R}^3\setminus \{O\} \to \mathbb S^2\times\mathbb S^1$ given by the formula $$\pi(x)=(s(x),r(x)\sim\,mod~1).$$  One can check directly that the standard frame $\mathbb F_\nu$ and $\mathbb R^3$ is projected by $\pi$ into the circles frame $$\mathbb J_\nu=\bigcup\limits_{j=1}^{\nu}\mathbb C_j,~~\mathbb C_j=\{q_j\}\times\mathbb S^1.$$

We say the circles frame $\mathbb J_\nu$ to be {\it standard}. Any collection $C_1,\dots, C_\nu$ of pairwise disjoint topologically embedded circles which are generators of the fundamental group $\pi_1(\mathbb S^2\times\mathbb S^1)$, forms a {\it frame of $\nu$ circles} $$J_\nu=\bigcup\limits_{j=1}^{\nu}C_j.$$ 

We say a frame $J_\nu$ to be {\it trivial} if there is a homeomorphism $\Phi:\mathbb{S}^2\times\mathbb S^1\to \mathbb{S}^2\times\mathbb S^1$ such that  ${J}_\nu=\Phi(\mathbb{J}_\nu)$. 

\begin{figure}[h]
\centerline{\includegraphics[width=6 true cm]{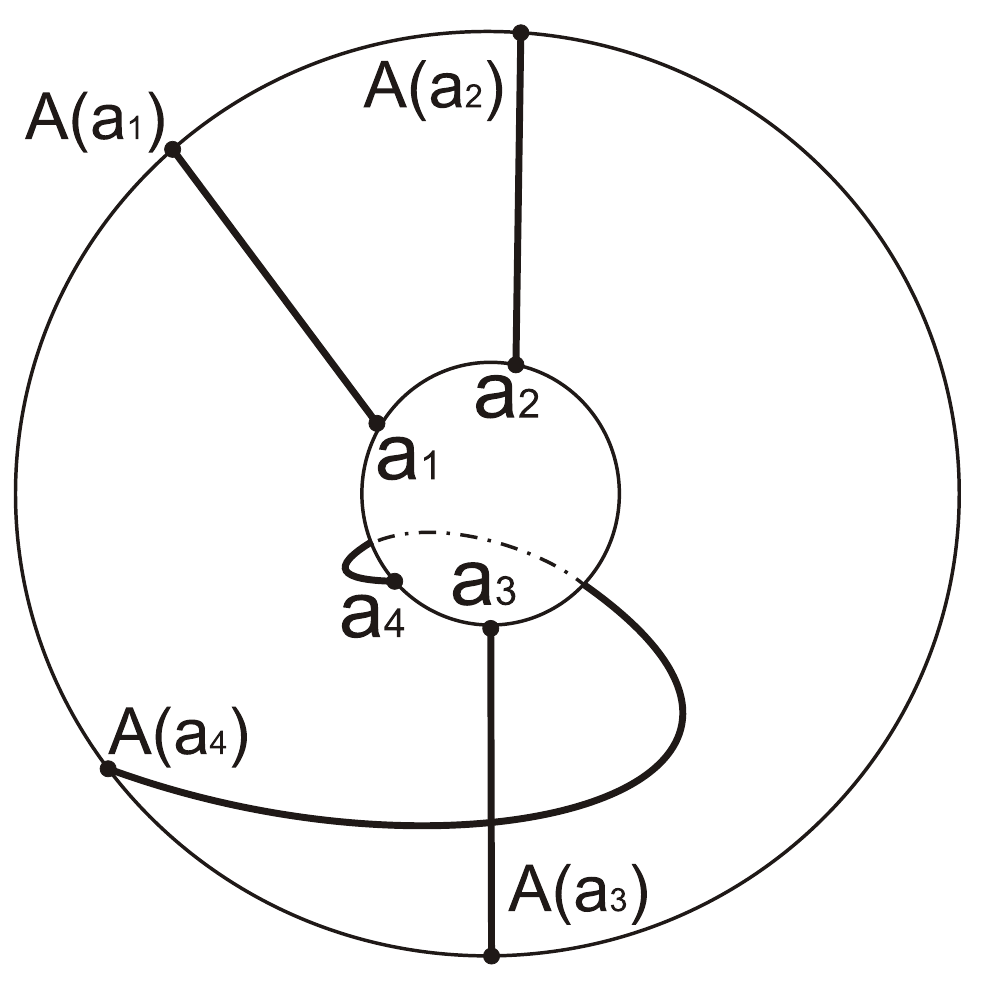}}
\caption{An example of a nontrivial frame of four circles that generates a tame frame of arcs}
\label{Sh+5}
\end{figure}

Every circles frame $J_\nu$ lifts to the arcs frame $F_{J_\nu}=\pi^{-1}(J_\nu)\cup O$ and if $J_\nu$ is trivial then according to Theorem \ref{fibr} the corresponding arcs frame is tame. The opposite generally is not true. In Lemma 5.1 of \cite{GGMP2012} a nontrivial circles frame $J_4$ was constructed for which the corresponding $F_{J_4}$ is tame (see Fig. \ref{Sh+5}). This construction can be generalized for any $\nu\geq 3$. On the other hand, for $\nu=1$ triviality of the circle is equivalent to tameness of the corresponding arc \cite[Lemmas 4.3, 4.4]{GrMedPoch2016}. In this paper we prove that the same holds for $\nu=2$.

\begin{theorem}\label{two} Triviality of a circles frame $J_2$ is equivalent to tameness of the corresponding arcs frame $F_{J_2}$. 
\end{theorem}

The key point in the proof of Theorem \ref{two} is the following lemma which is also a stand alone result.
\begin{lemma}\label{lemma} Let $J_\nu$ be a circles frame lifting to the tame arcs frame $F_{J_\nu}$. Then there exists a homeomorphism $G_1:\mathbb R^3\setminus O\to\mathbb R^3\setminus O$ such that $G_1A=AG_1$ and every sphere of the foliation $P=G_1(\mathbb P)$ intersects each arc $L_j,j\in\{1,\dots,\nu\}$ from $F_{J_\nu}$ at single point.
\end{lemma}

Now we use theorems \ref{fibr} and \ref{two} to investigate an embedding of frames of one-dimensional separatrices in basins of sources of rough 3-diffeomorhisms with two-dimensional expanding attractor.

In more details.  

Let $M^3$ be a closed 3-manifold and $f:M^3\to M^3$ be a structural stable diffeomorphism. Due to [Smale 1967] the non-wandering set $NW(f)$ of $f$ consists of a finite union of pairwise disjoint closed invariant sets  $\Lambda_1, \ldots, \Lambda_k$ called {\it basis sets}, each of them contains a dense orbit. A basic set is called {\it non-trivial}, if it is not a periodic orbit (in particular, not a fixed point).

A basic set $\Lambda$ is called {\it an attractor} if $\Lambda=\bigcup\limits_{x\in\Lambda}W^u(x)$. An attractor $\Lambda$ of the diffeomorphism  $f$ is called an {\it expanding attractor} if its topological dimension equals the dimension of the unstable manifold $W^u(x),x\in\Lambda$. 

Denote by $G$ a set of structural stable diffeomorphisms $f:M^3\to M^3$ with two-dimensional expanding attractor $\Lambda$. By [Grines V., Zhuzhoma E. On structurally stable diffeomorphisms with codimension one expanding attractors// Trans. Am. Math. Soc.—2005.— 357, 2.— P. 617-667.], other basic sets of $f\in G$ are trivial and its set $R$ of source points  is non empty. Moreover, with every periodic source $\alpha\in R$ is associated exactly two stable saddle separatrices $\ell_1^\alpha,\ell_2^\alpha$ which completely belong to the basin $W^u(\alpha)$, they have the same period, denote it $m_\alpha$. As $\alpha$ is a hyperbolic point then  the diffeomorphism $f^{m_\alpha}|_{W^u(\alpha)}$ is topologically conjugated with the homothety $A$ by means a homeomorphism $h_\alpha:W^u(\alpha)\to \mathbb R^3$. Let $L^\alpha_1=h_\alpha(\ell_1^\alpha),L^\alpha_2=h_\alpha(\ell_2^\alpha)$. Then the union  $$F^\alpha_2=L_1^\alpha\cup L_2^\alpha\cup O$$ is a  frame of two arcs, which we call {\it associated with the source $\alpha$}. Let $$J^\alpha_2=\pi(F^\alpha_2).$$ Then $J^\alpha_2$ is a  frame of two circles in $\mathbb S^2\times\mathbb S^1$, which we call {\it associated with the source $\alpha$}. Using theorems \ref{fibr} and \ref{two} we prove the following result.
   
\begin{theorem}\label{three} For every periodic source $\alpha\in R$ of a diffeomorphism $f\in G$ the associated with $\alpha$ frame of arcs $F^\alpha_2$ is tame and the frame of circles $J^\alpha_2$ is trivial. 
\end{theorem}

\section{Criterion for an arcs frame to be tame}

In this section we prove Theorem \ref{fibr}. Necessity follows from the definition of a tame frame. 

Sufficiency. Let there exist a homeomorphism $G:\mathbb R^3\setminus O\to\mathbb R^3\setminus O$ such that every sphere of the foliation $P=G(\mathbb P)$ intersects each arc $L_j,j\in\{1,\dots,\nu\}$ of the frame $F_\nu=\bigcup\limits_{j=1}^\nu L_j\cup O$ at a single point. Let $\tilde L_j=p(G^{-1}(L_j))$ and $\tilde p_{j,t}=\tilde L_j\cap(\mathbb S^2\times\{t\})$. Then $\tilde p_{j,t}$ has coordinates $(\tilde q_{j,t},t)$ in $\mathbb S^2\times\mathbb R$. 

Define the isotopy $\xi_{t}: \tilde q_{1,0}\cup\dots\cup\tilde q_{\nu,0}\to\mathbb S^2,t\in[0,1]$ by $$\xi_t(\tilde q_{j,0})=\tilde q_{j,t}.$$ By construction $\xi_0(\tilde q_{j,0})=\tilde q_{j,0}$. According to Thom's isotopy extension Theorem \cite[Theorem 5.8]{Mil1965} there is an isotopy $\chi_{t}:\mathbb S^2\to\mathbb S^2,t\in[0,1]$ coinciding with $\xi_{t}$ on $\tilde q_{1,0}\cup\dots\cup\tilde q_{\nu,0}$ and such that $\chi_0=id$. Let $q_1,\dots,q_\nu\subset\mathbb S^2$ be a collection of pairwise disjoint points. According to \cite[Ch. 8, Theorem 3.1]{hir} there is a diffeomorphism $\eta_0:\mathbb S^2\to\mathbb S^2$ which is isotopic to the identity and such that $\eta_0(q_j)=\tilde q_{j,0}$. Then the map $\eta_{[0,1]}:\mathbb S^2\times[0,1]\to\mathbb S^2\times[0,1]$ defined by  $$\eta_{[0,1]}(s,t)=(\chi_t(\eta_0(s)),t),$$ is a diffeomorphism and it sends the intervals $\{q_j\}\times[0,1]$ to the arcs $\tilde L_j\cap(\mathbb S^2\times[0,1])$. Define the diffeomorphism  $\eta_1:\mathbb S^2\to\mathbb S^2$ by $$\eta_1(s)=\chi_1(\eta_0(s)).$$ Then we similarly extend the diffeomorphism $\eta_1$ to all the sets $\mathbb S^2\times[i, i+1]$, $\mathbb S^2\times[-i,-i+1]$, $i\in \mathbb N$ one by one and this concludes the proof.

\section{Proof of Lemma \ref{lemma}}

Let $J_\nu$ be a circles frame which lifts to a frame of tame arcs $F_{J_\nu}$. In this section we show that there is a homeomorphism $G_1:\mathbb R^3\setminus O\to\mathbb R^3\setminus O$ such that $G_1A=AG_1$ and every sphere of the foliation $P_1=G_1(\mathbb P)$ intersects each arc $L_j,j\in\{1,\dots,\nu\}$ of the frame $F_{J_\nu}$ at a single point.

\begin{proof} Let $S_0=G_1(\mathbb S^2)$. It follows from the definition of a tame arcs frame that $S_0$ is the 2-sphere tamely embedded into $\mathbb R^3$ and it 

intersects each arc $L_j,j\in\{1,\dots,\nu\}$ at a single point\hfill (*) 

According to \cite{hir} we can think that the sphere $S_0$ is smooth also as any compact parts of arcs $L_j$.  Denote by $B_0$ the smooth 3-ball in $\mathbb R^3$ bounded by $S_0$. The ball $B_0$ will now be transformed into a ball $B$ such that $A^{-1}(B)\subset int\,B$ and the sphere $S=\partial B$ satisfies (*).
 
Let $m$ be the least natural for which $A^{-k}(S_0)\cap S_0=\emptyset$ for each $k>m$. If $S_0$ is not transversal to its images $A^{-1}(S_0),\dots,A^{-m}(S_0)$ we modify it in the following way. For every $x\in S_0$ choose a compact neighborhood $K_x$ of the point $x$ in $\mathbb R^3$ such that $A^{-1}(K_x)\cap K_x=\emptyset$. Such a neighborhood exists because there are no fixed points of the diffeomorphism $A^{-1}$ on the sphere $S_0$. Pick a finite subcover $K_1,\dots,K_p$ of the sphere $S_0$ from a cover $\{K_x,~x\in S_0\}$. According to Transversality Theorem \cite{hir} $S_0$ can be approximated by such a smooth sphere (we denote it by $S_0$ again) that $A^{-1}(K_1)$ is transversal to $S_0$. Then we do the same for  $K_1\cup K_2$ ans so on up to $K_1\cup\dots\cup K_p$. In the same way we approximate the sphere so that it is transversal to its images by $A^{-1}$ and $A^{-2}$. After $m$ steps we get a sphere (denoted again by $S_0$) which is transversal to all its images by the maps $A^{-1}, \dots, A^{-m}$. All these approximations can be so chosen that $S_0$ satisfies  (*). 

Now $S_0$ is transversal to its images $A^{-1}(S_0),\dots,A^{-m}(S_0)$. There are two cases: 
\begin{enumerate}[(I)]
\item $m=1$, i.e. $A^{-1}(S_0)\cap S_0\neq\emptyset$, $A^{-k}(S_0)\cap S_0=\emptyset$ for $k>1$,
\item $m>1$, i.e. $A^{-k}(S_0)\cap S_0\neq\emptyset$ for $k=1,\dots,m$;\ $A^{-k}(S_0)\cap S_0=\emptyset$ for $k>m$.
\end{enumerate}

Consider case (I). The sphere $S_0$ will now be modified into a sphere $S'_0$ satisfying (*) and the following two conditions:

\begin{enumerate}[(i)]

\item $A^{-k}(S'_0)\cap S'_0=\emptyset$ for $k>1$;

\item the number of the arcs in the intersection $A^{-1}(S'_0)\cap S'_0$ is less than the number of the arcs in the intersection $A^{-1}(S_0)\cap S_0$. 

\end{enumerate}

Let  $\Sigma=A^{-1}(S_0)$, $\tilde S_0=A^{-2}(S_0)$ and $\tilde B_0=A^{-2}(B_0)$. Let $\gamma$ be an arc in the intersection $\Sigma\cap(S_0\cup\tilde S_0)$. We say the arc $\gamma$ to be {\it extreme} on $\Sigma$ if it bounds such a disk $D_\gamma\subset\Sigma$ that $int~D_\gamma$ contains no arcs from $\Sigma\cap(S_0\cup\tilde S_0)$. Let $\gamma$ be this extreme arc. There are two cases: a) $\gamma\subset S_0$ and b) $\gamma\subset\tilde S_0$. 

In the case a) consider two subcases: a1) $D_\gamma\cap int~B_0=\emptyset$ and a2) $D_\gamma\subset B_0$. 

In the subcase a1) the disk $D_\gamma$ divides the domain $\mathbb R^3\setminus int~B_0$ into two parts and the closure in $\mathbb R^3$ of one of these parts is the 3-ball (denote it by $B_\gamma$). By construction $int~ D_\gamma\cap A^{-k}(S_0)=\emptyset$ for any $k>1$ and the disks $D_\gamma$ and $B_\gamma\cap S_0$ contain the same number (0 or 1) of points of intersection with the arc $L_j$. Therefore, by smoothing the angles of $B_0\cup B_\gamma$ one gets a 3-ball $B'_0$ whose boundary is the desired sphere $S'_0$ (see Fig. \ref{Bgamma}).

\begin{figure}[h]\label{Bgamma}
\centerline{\includegraphics[width=7 true cm]{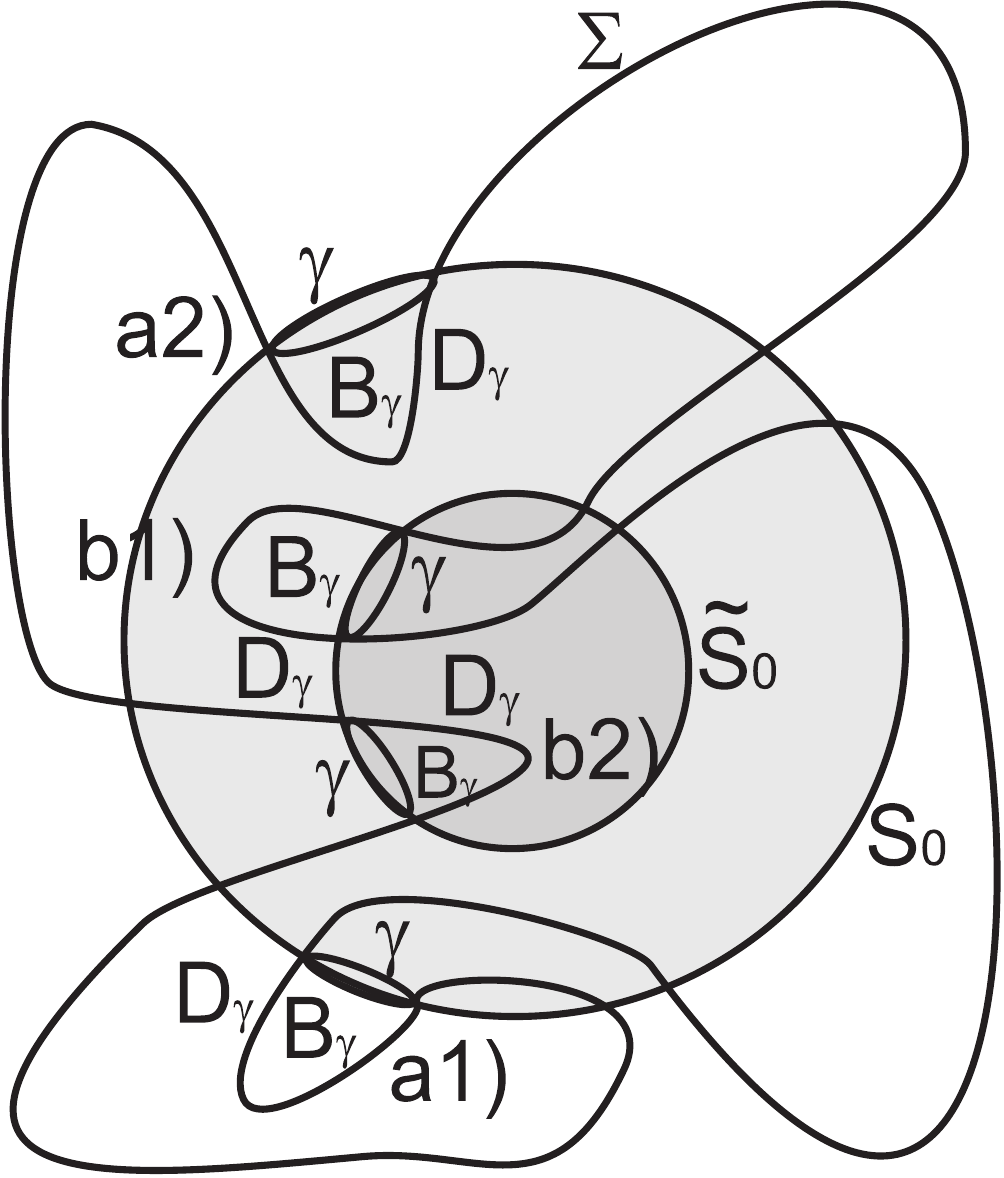}}
\caption{Illustration to the proof of Lemma \ref{lemma}} 
\end{figure}

In the subcase a2) the disk  $D_\gamma$ divides the 3-ball $B_0$ into two 3-balls. Only one of these balls does not contain $O$, because $O\notin\Sigma$ (denote its closure by $B_\gamma$). Similarly to the subcase a1) one gets the desired sphere as the boundary of the 3-ball $B'_0$ after smoothing the angles of $cl(B_0\setminus B_\gamma)$.  

In the case b) we again have two subcases: b1) $\tilde D_\gamma\cap int~\tilde B_0=\emptyset$ and b2) $\tilde D_\gamma\subset\tilde B_0$. 

In the subcase b1) the disk $\tilde D_\gamma$ divides the domain $\mathbb R^3\setminus int~\tilde B_0$ into two parts. The closure in $\mathbb R^3$ of one of these parts is the 3-ball (denote it by $\tilde B_\gamma$). By construction $int~\tilde D_\gamma\cap A^{-k}(\tilde S_0)=\emptyset$ for any $k<-1$ and the disks $\tilde D_\gamma$ and $\tilde B_\gamma\cap\tilde  S_0$ contain the same number (0 or 1) of points of intersection with the arc $L_j$. Therefore, by smoothing the angles of the ball $\tilde B_0\cup \tilde B_\gamma$ one can get such a 3-ball $\tilde B'_0$ that its boundary $\tilde S'_0$ satisfies  (*) and the following two conditions:
\begin{enumerate}[($\tilde{\rm i})$]
\item $A^{-k}(\tilde S'_0)\cap\tilde S'_0=\emptyset$ for $k<-1$;
\end{enumerate}
\begin{enumerate}[($\tilde{\rm ii})$]
\item  the number of the arcs in the intersection $A^{-1}(\tilde S'_0)\cap\tilde S'_0$ is less than the number of the arcs in the intersection $A^{-1}(\tilde S_0)\cap\tilde S_0$.
\end{enumerate}
Then the boundary of the ball $A^{-2}(\tilde B'_0)$ is the desired sphere $S'_0$.

In the subcase b2) the disk  $\tilde D_\gamma$ divides the 3-ball $\tilde B_0$ into two 3-balls. Exactly one of these balls does not contain the coordinates origin because $O\notin\Sigma$ (denote its closure by $\tilde B_\gamma$). Similarly to the subcase b1) the desired sphere is the boundary of the 3-ball $A^{-2}(\tilde B'_0)$ where the ball $\tilde B'_0$ is the result of smoothing the angles of the ball $cl(\tilde B_0\setminus \tilde B_\gamma)$. 

If we continue this procedure we get the smooth sphere $S$ which meets the arc $L_j$ at a single point and which bounds the ball $B\ni O$ such that $A^{-1}(B)\subset int~ B$.

Consider the case (II). Pick a natural $r$ such that $2^r\leq m<2^{r+1}$. Let $g_r=A^{-2^r}$ then $g_r^k(S_0)\cap S_0=\emptyset$ for all $k\geq 2$. Similarly to the previous case we construct a sphere $S_1$ which intersects the arc  $L_j$ at a single point and  $g_r(S_1)\cap S_1=\emptyset$. Thus we have decreased $r$ at least by 1: $A^{-2^r}(S_1)\cap S_1=\emptyset$. We continue this procedure and we get the desired ball $B$.

Now we are going to show that the annulus $K=B\setminus int\,A^{-1}(B)$ can be foliated into 2-spheres each of whish satisfies (*).

Choose an annulus  $\tilde K=G_1(\mathbb S^2\times[-n,n]),n\in\mathbb N$ containing the annulus $K$. Let $S_t=G_1(\mathbb S^2\times\{t\},t\in[-n,n]$. Choose tubular neighborhoods $U_j$ of the arcs $L_j$ such that each of $U_j\cap S_t$, $U_j\cap S$ and  $U_j\cap A^{-1}(S)$ consists of one 2-disk. Thus, the manifold $M=\tilde K\setminus(\bigcup\limits_{j=1}^\nu U_j)$ is homeomorphic to $\Sigma_\nu\times[-n,n]$ where $\Sigma_\nu\times\{t\}=S_t\setminus(\bigcup\limits_{j=1}^\nu U_j)$ is the 2-sphere with $\nu$ holes.  Let $\Sigma^1=S\setminus(\bigcup\limits_{j=1}^\nu U_j)$ and  $\Sigma^2=A^{-1}(S)\setminus(\bigcup\limits_{j=1}^\nu U_j)$. The surface $\Sigma^1$ is the 2-sphere  with $\nu$ holes as well and it separates in  $M$ the surfaces $\Sigma_\nu\times\{-n\}$ and $\Sigma_\nu\times\{n\}$. According to \cite{GrMeZh2003} the surface $\Sigma_1$ divides the manifold $M$ into two parts, one of which is homeomorphic to $M$. By the same reasons the surfaces  $\Sigma^1$ and $\Sigma^2$ bound in  $M$ a manifold homeomorphic to $\Sigma_\nu\times[0,1]$. Therefore the annulus $K$ has the desired properties. 
\end{proof}

\section{Criterion for a frame of two circles to be trivial}

In this section we prove Theorem \ref{two}. To do this it is sufficient to show that if an arcs frame $F_{J_2}$ is tame then the circles frame $J_2$ is trivial. By Theorem \ref{fibr} if an arcs frame $F_{J_2}=L_1\cup L_2\cup O$ is tame then there is a homeomorphism $G_1:\mathbb R^3\setminus O\to\mathbb R^3\setminus O$ such that every sphere of the foliation  $P_1=G_1(\mathbb P)$ intersects each arc $L_j,j\in\{1,2\}$ at a single point. 

Let $\tilde L_j=G_1^{-1}(L_j)$. Let $\mathbb L_j$ be a standard ray passing through $q_j=\tilde L_j\cap\mathbb S^2$. Let $K=p^{-1}(\mathbb S^2\times[0,1])$ and let $\Pi$ be the plane in $\mathbb R^3$ containing the rays $\mathbb L_1$ and $\mathbb L_2$. Denote by $H\subset\Pi$ one of the sectors bounded by these rays on $\Pi$ and let $I_j=\mathbb L_j\cap K,\tilde I_j=\tilde L_j\cap K$. Let $V(I_i)\subset K$ be a tubular neighborhood of the segment $I_i$ whose intersections with every sphere $p^{-1}(\mathbb S^2\times\{t\})$ is a disk and let $V(\tilde I_i)$ be a similar neighborhood of $\tilde I_j$ and such that  $V(I_j)\cap\partial K=V(\tilde I_j)\cap\partial K$. Let $W=K\setminus(V(I_1)\cup V(I_2)),\tilde W=K\setminus(V(\tilde I_1)\cup V(\tilde I_2))$, $D=H\cap W,I^j=\partial D\cap V(I_j), E=D\cap \mathbb S^2$. By construction $W$ is the solid torus, $\mu=E\cup I^1\cup A(E)\cup I^2$ being its meridian. The manifold $\tilde W$ is the solid torus as well. Define a homeomorphism $\psi:\partial W\to\partial{\tilde W}$ in such a way that 
\begin{enumerate}
\item it coincides with the identity on $\partial K\cap\partial W$,
\item it coincides on $\partial V(I_j)$ with the homeomorphism that sends disks to disks and such that $\psi(\mu)$ is  a meridian of the solid torus $\tilde W$. 
\end{enumerate}
By \cite[p. 30]{Ro} the homeomorphism $\psi$ can be extended to a homeomorphism $\tilde\Psi:K\to K$ commuting with $A$. Therefore it projects into the desired homeomorphism $\Psi:\mathbb{S}^2\times\mathbb S^1\to \mathbb{S}^2\times\mathbb S^1$ such that  $\Psi(\pi(\mathbb{L}_1\cup\mathbb L_2))=\pi(\tilde L_1\cup\tilde L_2)$.

\section{Application to dynamics. Proof of the theorem \ref{three}}

Denote by $G$ a set of structural stable diffeomorphisms $f:M^3\to M^3$ whose non-wandering set $NW(f)$ contains a two-dimensional expanding attractor $\Lambda$. Firstly, we describe a dynamics of a diffeomorphism $f\in G$ following by [Grines V., Zhuzhoma E. On structurally stable diffeomorphisms with codimension one expanding attractors// Trans. Am. Math. Soc.—2005.— 357, 2.— P. 617-667.].

For any point $x\in \Lambda$, $\dim W^s(x)=1$ allows one to introduce  the notation  $(y,z)^s$ ($[y,z]^s$) for an open (closed) arc of the  stable manifold $W^s(x)$ bounded by points $y,z\in W^s(x)$. 
The set $W^s(x)\backslash x$ consists of two components. At least one of these components has a nonempty intersection with the set $\Lambda$. A point $x\in \Lambda$ is called {\it boundary} if one of the connected components of the set $W^s(x)\backslash x$ does not intersect with $\Lambda$, we will denote this component by $W^{s\emptyset}(x)$. 

The set $\Gamma_{\Lambda}$ of all boundary points  of the set $\Lambda$ is nonempty  and consists of finite number of periodic points  that are divided into associated couples $(p,q)$ of points of the same period so the 2-bunch $B_{pq}=W^u(p)\cup W^u(q)$ is a boundary achievable\footnote{For an open set $D\subset M^3$ with boundary $\partial D$ a subset $\delta D\subset \partial D$ is called {\it achievable from  inside of the domain $D$} if for any point $\delta\in \delta D$ there exists an arc $l_{\delta, d}$ with endpoints $\delta$ and $d\in D$  such that  $(l_{\delta, d}\setminus\delta)\subset D$.} from inside of the connect component of the set $M\backslash \Lambda$. For each couple $(p,q)$  of associated boundary points of the set $\Lambda$ we construct the so-called characteristic sphere.
\begin{figure}[h]
\begin{center}
\includegraphics[width=10cm,height=6cm]{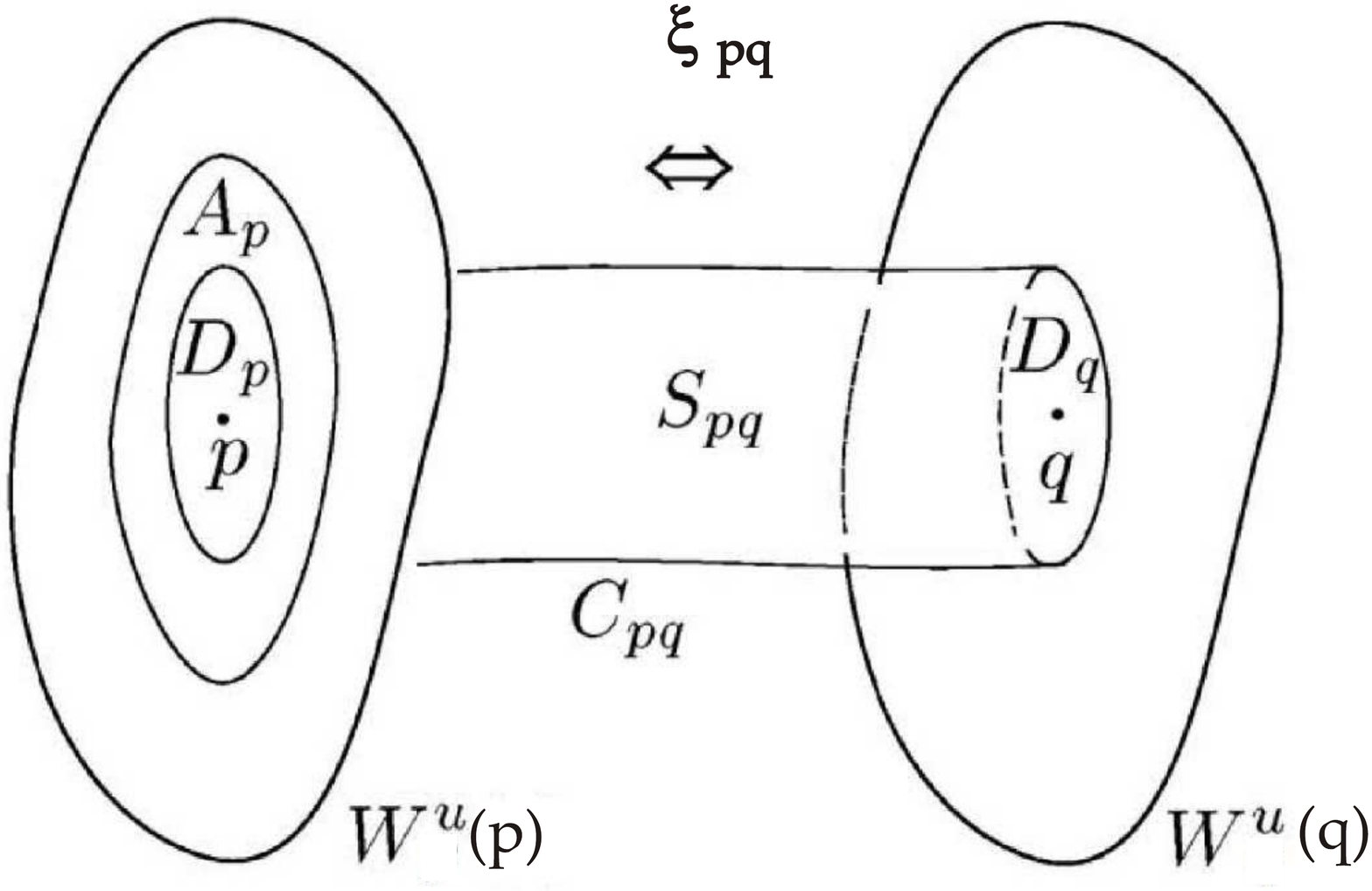}
\end{center}
\caption{Characteristic sphere}
\label{sphere}
\end{figure}

Let $B_{pq}$ be a 2-bunch of  the attractor $\Lambda$ consisting of two unstable manifolds $W^u(p)$ and $W^u(q)$ of associated boundary points  $p$ and $q$ respectively and let $m_{pq}$  be the period of points $p$, $q$. Then for any point $x\in W^u(p) \backslash p$ there exists a unique point $y\in (W^u(q)\cap W^s(x))$ such that the arc $(x,y)^s$ does not intersect with the set $\Lambda$. We define the mapping 
$$\xi _{pq} : B_{pq} \backslash \{p,q\}\rightarrow B_{pq} \backslash  \{p,q\} $$ 
by putting $\xi _{pq}(x)=y$ and $\xi _{pq}(y)=x$. Then $\xi _{pq}(W^u(p)\backslash p)=W^u(q)\backslash q$, i.e. the mapping  $\xi _{pq}$ takes the pierced unstable manifolds of the 2-bunch into each other and is an
 involution ($\xi^2 _{pq}(x)=id$). By the theorem on continuous dependence of invariant manifold on compact sets the mapping $\xi _{pq}$  is a homeomorphisms. 
 
The restriction $f^{m_{pq}}|_{W^u(p)}$ has exactly one hyperbolic repelling fixed point $p$, hence there exists a smooth 2-disk $D_p\subset W^u(p)$ such that $p\in f^{-m_{pq}}(D_p)\subset int\,D_p$. Then  the set $C_{pq}=\bigcup\limits_{x\in \partial D_p}(x,\xi _{pq}(x))^s $ is homeomorphic to a closed cylinder $\mathbb{S}^1\times[0,1]$. The set $C_{pq}$ is called a {\it connecting  cylinder}. The circle $\xi _{pq}(\partial D_p)$ bounds in $W^u(q)$ a two-dimensional 2-disk  $D_q$ such that $q\in D_q\subset  int(f^{m_{pq}}(D_q))$. The set $S_{pq}=D_p\cup C_{pq}\cup D_q$ is homeomorphic to a 2-dimensional sphere called the {\it characteristic  sphere} corresponding to the bunch $B_{pq}$  (fig. \ref{sphere}).
 
\begin{figure}[h]
\begin{center}
\includegraphics [width=11cm,height=7cm]{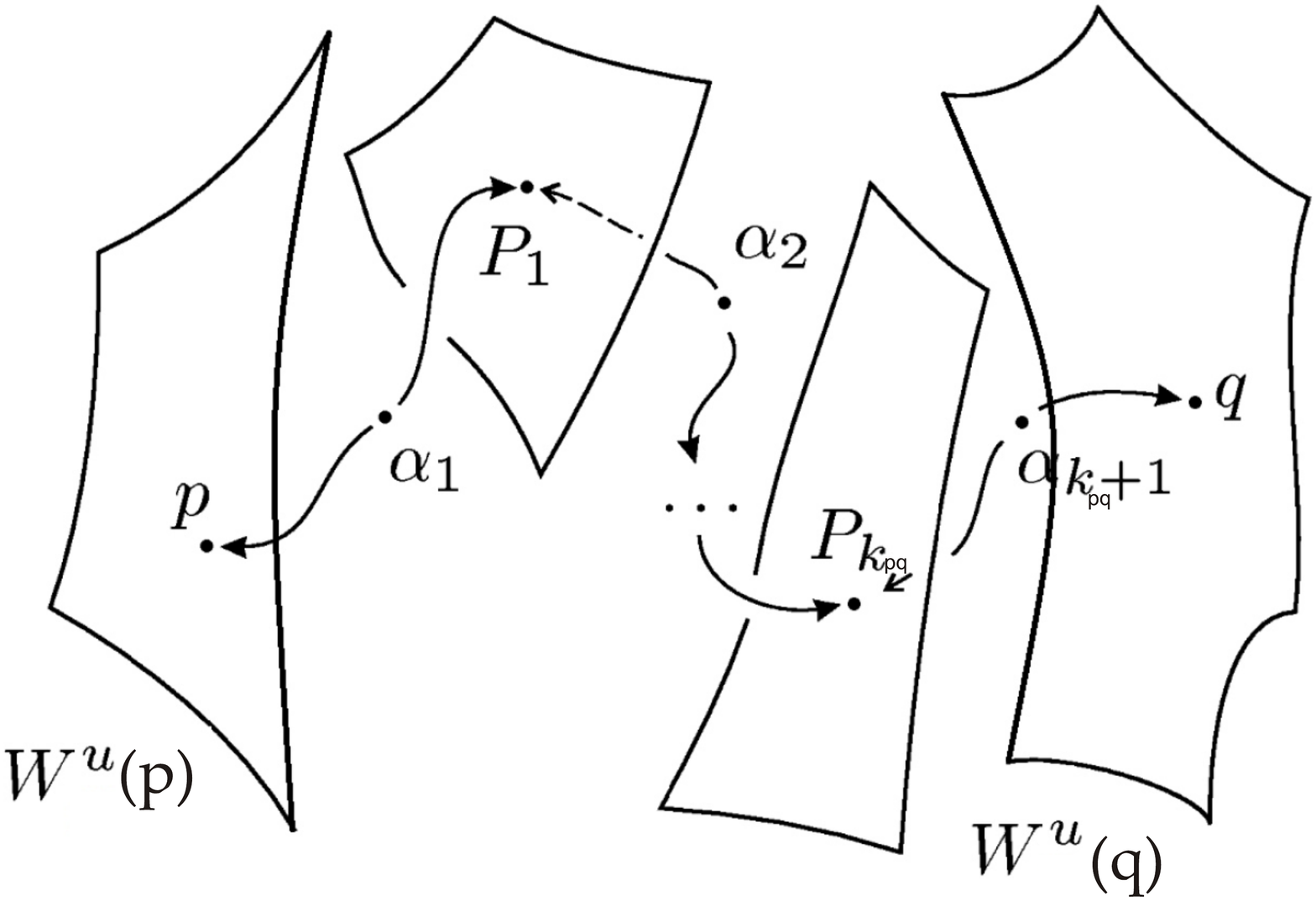}
\end{center}
\caption{Arc $l_{pq}$}
\label{curve}
\end{figure}

\begin{proposition}[Grines V., Zhuzhoma E. On structurally stable diffeomorphisms with codimension one expanding attractors// Trans. Am. Math. Soc.—2005.— 357, 2.— P. 617-667.] For every diffeomorphism $f\in G$ with two-dimensional expanding attractor $\Lambda$ the following facts take place:

(1) the ambient manifold $M^3$ is homeomorphic to the 3-dimensional torus $\mathbb{T}^3$;

(2) each characteristic sphere $S_{pq}$ bounds a 3-ball $Q_{pq}$ such that $(NW(f)\setminus\Lambda)\subset \bigcup \limits_{(p,q)\subset\Gamma_{\Lambda}}Q_{pq}$;

(3) for each associated couple $(p,q)$ of boundary points there exists a natural number $k_{pq}$ such that $NW(f)\cap Q_{pq}$ consists of $k_{pq}$ periodic sources $\alpha_1,...,\alpha_{k_{pq}}$ and $k_{pq}-1$ periodic saddle points $P_1,...,P_{k_{pq}-1}$ alternate  on the simple arc $l_{pq}=W^{s\emptyset }(p)\cup \bigcup\limits_{i=1}^{k_{pq}-1}W^s(P_i)\cup \bigcup\limits_{i=1}^{k_{pq}}W^s(\alpha_i)\cup W^{s\emptyset}(q)$ (fig. \ref{curve}).
\end{proposition}  
 
Let as assume for completeness $P_0=p$ and $P_{k_{pq}}=q$. For $i\in\{1,\dots,k_{pq}\}$ let $$\alpha=\alpha_i,\,p_1=P_{i-1},\,p_2=P_i.$$  
Due to structural stability of diffeomorphism $f$ each arc $(x,y)$ where $x\in (D_p\setminus p)$, $y\in (D_q\setminus q)$ intersects $W^u(p_j),j=1,2$ exactly at one point. Indeed, assuming the converse, we find a point $z\in\Lambda$ such that $W^s(z)$ touches $W^u(p_j)$, which contradicts the strong transversality condition. Then 3-dimensional disk $Q_{pq}$ intersects the two-dimensional unstable manifold of the saddle $p_j$ exactly at one two-dimensional disk $D_{p_j}$. 

\begin{figure}[h]
\begin{center}
\includegraphics [width=13cm,height=7cm]{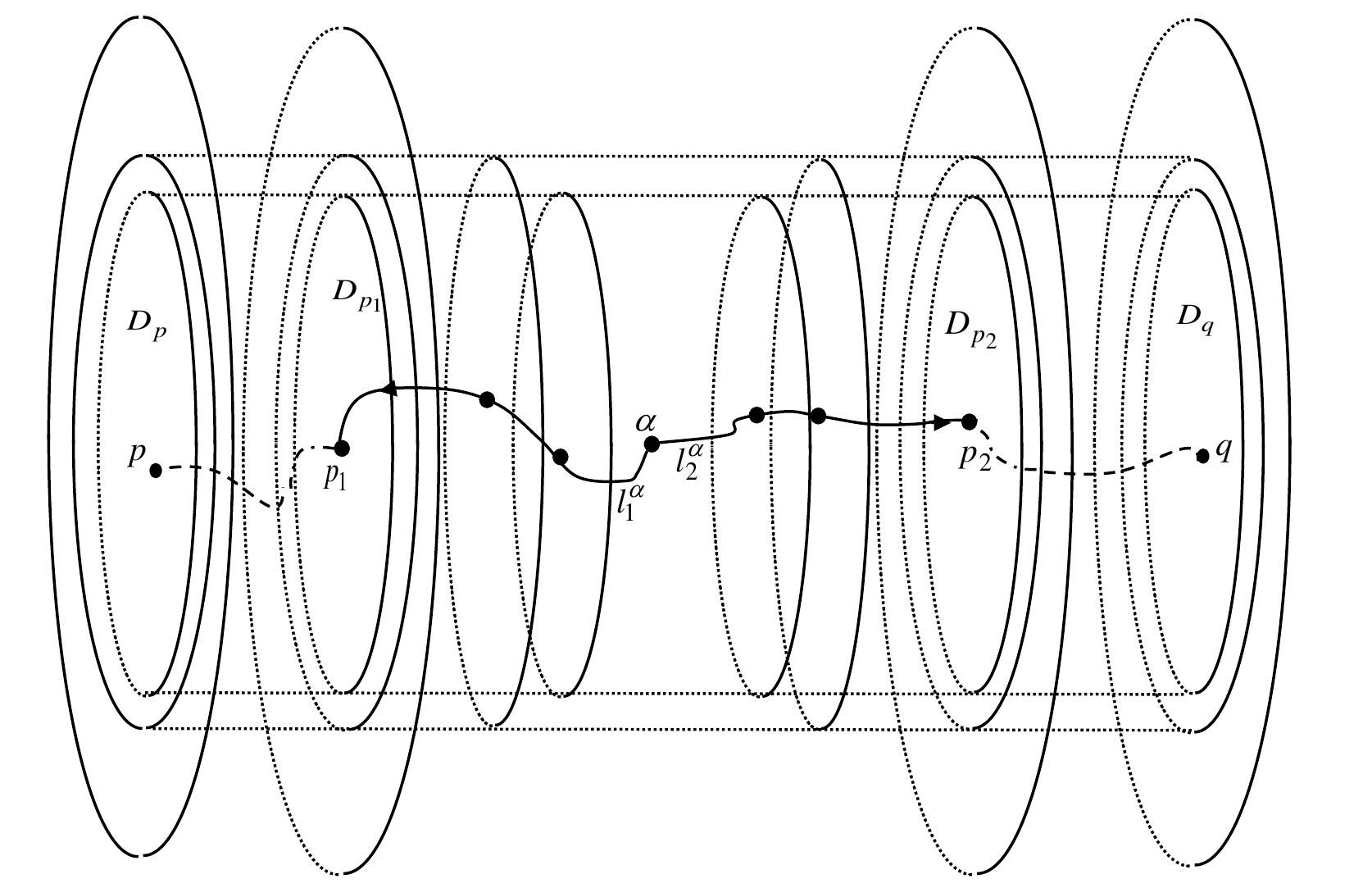}
\end{center}
\caption{Illustration to the proof of the theorem \ref{three}}
\label{annulus}
\end{figure}

Denote by $Q_\alpha$ the closure of the connected component of the set $Q_{pq}\setminus(D_{p_1}\cup D_{p_2})$ containing the source $\alpha$. Let $S_\alpha=\partial Q_\alpha$ and $\ell^\alpha_j=W^s(p_j)\cap W^u(\alpha)$  (fig. \ref{annulus}). Notice that the period $m_\alpha$ of the separatrix $\ell^\alpha_j$ equals $m_{pq}$. 

As $\alpha$ is a hyperbolic point then  the diffeomorphism $f^{m_\alpha}|_{W^u_\alpha}$ is topologically conjugated with the homothety $A$ by means a homeomorphism $h_\alpha:W^u_\alpha\to \mathbb R^3$. Let $L^\alpha_1=h_\alpha(\ell_1^\alpha),L^\alpha_2=h_\alpha(\ell_2^\alpha)$. Then the union  $F^\alpha_2=L_1^\alpha\cup L_2^\alpha\cup O$ is associated with the source $\alpha$ frame of two arcs. Let $J^\alpha_2=\pi(F^\alpha_2)$, then $J^\alpha_2$ is associated with the source $\alpha$ frame of two circles in $\mathbb S^2\times\mathbb S^1$. Due to theorems \ref{fibr} and \ref{two} to prove that the associated with $\alpha$ frame of arcs $F^\alpha_2$ is tame and the frame of circles $J^\alpha_2$ is trivial it is enough 
to construct a foliation on the manifold $W^u(\alpha)\setminus\alpha$ each leaf of which is a 2-spheres around $\alpha$ intersecting every separatrices $\ell^\alpha_1$ and $\ell^\alpha_2$ exactly once. 

For this aim let us consider a point $x_j\in\ell^\alpha_j$ and a fundamental segment $I_j\subset\ell^\alpha_j$ with endpoint $x_j,f^{m_\alpha}(x_j)$. Let us choose a tubular neighborhood $V(I_j)$ which is diffeomorphic to $\mathbb D^2\times[0,1]$ by means a diffeomorphism $h_j:\mathbb D^2\times[0,1]\to V(I_j)$ such that $f^{m_\alpha}(h_j(\mathbb D^2\times\{0\}))\subset h_j(\mathbb D^2\times\{1\})$. Due to $\lambda$-lemma we can assume that every disc $h_j(\mathbb D^2\times\{0\})$ intersects $W^s(x)$ for every $x\in (D_p\setminus p)$ exactly once (in the opposite case we can choose an appropriate iteration of $V(I_j)$). 

The set $K_p=D_p\setminus int\,f^{-m_\alpha}(D_p)$ is diffeomorphic to $\mathbb S^1 \times [0,1]$ by means a diffeomorphism $h:\mathbb S^1 \times [0,1]\to K_p$. Denote by $S_t,t\in[0,1]$ a two-dimensional sphere bounded by the discs $h_1(\mathbb D^2\times\{t\}),\,h_2(\mathbb D^2\times\{t\})$ and the cylinder $h(\mathbb S^2 \times\{t\})$. Iterations of these spheres with respect to $f^{km_\alpha},k\in\mathbb Z$ give the required foliation on $W^u(\alpha)\setminus\alpha$.

\end{document}